\documentclass[12pt]{article}
\usepackage{titlesec}
\usepackage{lineno}
\usepackage{lipsum}
\usepackage{changepage}
\usepackage{graphicx}
\usepackage{amssymb}
\usepackage{amsmath}
\usepackage{amsfonts}
\usepackage{multirow}
\usepackage{booktabs}
\usepackage{cite}
\usepackage{caption}
\usepackage{setspace}
\usepackage{color}
\usepackage{comment}
\usepackage{epstopdf}
\usepackage{authblk}
\usepackage{bm}
\allowdisplaybreaks[4]
\usepackage{textcomp}
\newtheorem{thm}{\textbf{Theorem}}
\newtheorem{assu}{\textbf{Assumption}}
\newtheorem{rem}{\textbf{Remark}}
\newtheorem{lemma}{\textbf{Lemma}}
\newtheorem{corollary}{\textbf{Corollary}}
\usepackage[numbers,sort&compress]{natbib}

\usepackage{times}

\usepackage[
 left=3cm,  
 right=2.6cm,  
 top=2.1cm,   
 bottom=2.1cm, 
]{geometry}

\definecolor{qired}{rgb}{0.6, 0, 0}

\DeclareCaptionFormat{boldnumberfirstsentence}{
 \textbf{#1#2}#3\par
}

\DeclareCaptionLabelSeparator{colonspace}{: }
\newcommand*\patchAmsMathEnvironmentForLineno[1]{%
 \expandafter\let\csname old#1\expandafter\endcsname\csname #1\endcsname
 \expandafter\let\csname oldend#1\expandafter\endcsname\csname end#1\endcsname
 \renewenvironment{#1}%
   {\linenomath\csname old#1\endcsname}%
   {\csname oldend#1\endcsname\endlinenomath}}%
\newcommand*\patchBothAmsMathEnvironmentsForLineno[1]{%
 \patchAmsMathEnvironmentForLineno{#1}%
 \patchAmsMathEnvironmentForLineno{#1*}}%
\AtBeginDocument{%
\patchBothAmsMathEnvironmentsForLineno{equation}%
\patchBothAmsMathEnvironmentsForLineno{align}%
\patchBothAmsMathEnvironmentsForLineno{flalign}%
\patchBothAmsMathEnvironmentsForLineno{alignat}%
\patchBothAmsMathEnvironmentsForLineno{gather}%
\patchBothAmsMathEnvironmentsForLineno{multline}%
}

\allowdisplaybreaks

\title{High Probability Convergence of Distributed Clipped Stochastic Gradient Descent with Heavy-tailed Noise}
\date{}

\author[1]{\fontsize{12}{14}\selectfont Yuchen Yang}
\author[2]{\fontsize{12}{14}\selectfont Kaihong Lu}
\author[1,*]{\fontsize{12}{14}\selectfont Long Wang}
\affil[1]{Center for Systems and Control, College of Engineering, Peking University, Beijing, 100871, China.}
\affil[2]{College of Electrical Engineering and Automation, Shandong University of Science and Technology,
 Qingdao 266590, China}
\affil[*]{Corresponding author. Email: \text{longwang@pku.edu.cn}}

\begin{document}
\maketitle

\begin{abstract}
In this paper, the problem of distributed optimization  is studied via
 a network of agents. Each agent only has access to a noisy gradient of its own objective function,
 and can communicate with its neighbors via a network. To handle this problem, a distributed
 clipped stochastic gradient descent algorithm is proposed, and  the  high probability convergence of the algorithm is studied. Existing works on distributed algorithms involving stochastic gradients only consider the light-tailed
noises. Different from them, we study the case with heavy-tailed settings. Under mild assumptions on
 the graph connectivity, we prove that the algorithm converges in high probability under a certain clipping operator. Finally,
 a simulation is provided to demonstrate the effectiveness of our theoretical results.
\end{abstract}

% \newpage
\section{Introduction}
In distributed optimization, the goal of agents is to cooperatively minimize the global objective function formed
	by the sum of local functions \cite{4749425}. Along with the penetration of multi-agent networks \cite{wang2007new,wang2009controllability,5404832}, distributed optimization
	has received ever-increasing attention in recent years \cite{4749425,6930814,7518617}. This is due to its wide practical applications in many
	areas such as large-scale sensor networks \cite{7973152}, social networks \cite{2842722}, and distributed power systems \cite{9141512}.

 In practical applications, optimization problems usually
 occur in uncertain environments, and the accurate gradient is hard to obtain due to
uncertainties in communications and
environments. For instance, when using sensor networks to monitor a remote objective, e.g., avalanches \cite{4497285}, uncertainty is inevitable due to the measurement noise \cite{nedic2018network}. Thus, it is necessary to study the
distributed optimization problem with stochastic gradients. Recently, various stochastic gradient algorithms are proposed \cite{9241497,9084351,9910399,9669076}.
	In \cite{9910399}, a gradient tracking-based distributed optimization approach is proposed that can avoid information-sharing noise from accumulating in the gradient estimation.
	In \cite{9669076}, a  stochastic mirror descent algorithm with the event-triggered communication strategies is proposed. 

  It is worth noting that all the aforementioned investigations
 study the convergence of distributed stochastic algorithms
 in expectation. For the algorithm converging in expectation, it is necessary to run an algorithm in large numbers of rounds to eliminate the error between the unbiased estimator and the accurate gradient. However, in some practical problems, we can only run the algorithm for a few rounds  or even once, so achieving the high probability convergence is desired. For instance, in the problem of tracking a target, agents need to track the target as fast as possible. In fact, in centralized optimization, the high probability convergence has already been studied \cite{madden2021highprobability,hong2023high,pmlr-v202-liu23aa}. Unfortunately, the  results in \cite{madden2021highprobability,hong2023high,pmlr-v202-liu23aa} are not applicable to the distributed cases. More recently, high probability convergence is studied in distributed settings. In \cite{10295561}, the high probability convergence of a distributed stochastic gradient algorithm are provided. In \cite{10103556}, linear regression from data distributed over a network of agents by means of least absolute shrinkage and selection operator (LASSO) estimation is studied,  and the Distributed Gradient Descent in the Adapt-Then-Combine form is provided.
 
 All works \cite{madden2021highprobability,hong2023high,pmlr-v202-liu23aa,10295561,10103556} above focus on the light-tailed noise. 
For light-tailed noise, such as the sub-Gaussian noise, large variety of concentration
 techniques are applicable since its tail decays faster than exponential distribution \cite{10295561}. Even though the light-tailed noise assumption is intuitive, in domains like evolutionary ecology \cite{Jourdain2012LvyFI}, audio source separation \cite{7177973}, finance \cite{1997Fractals} and machine learning \cite{pmlr-v97-simsekli19a}, the assumption seems invalid, which implies the necessity of studying heavy-tailed setting.
 
\textcolor{black}{
To handle the heavy-tailed noise, the gradient clipping  strategy is employed, which is initially adopted in \cite{pmlr-v28-pascanu13} for the purpose of training neural networks. In Table \ref{table_comparison}, a brief review of the existing works and comparisons with our work are provided. For strongly convex objective functions, a centralized stochastic gradient algorithm is proposed in \cite{doi:10.1137/110848864}, where the convergence rate $\mathcal{O}(T^{-\frac{1}{2}})$ is achieved, while the distributed case is studied in \cite{9683110}. In \cite{doi:10.1137/110848864} and \cite{9683110}, the variance of the noises is assumed to be bounded. 
% In the case of finite variance,  and \cite{9683110} have separately demonstrated that, for centralized and distributed optimization problems, stochastic gradient descent can achieve a convergence rate of $\mathcal{O}(T^{-\frac{1}{2}})$ under the assumption of strong convexity of the objective function. 
In \cite{pmlr-v235-gorbunov24a,pmlr-v202-sadiev23a,nguyen2023high}, using constant step sizes based on the iteration number, the convergence rate of the clipped SGD is developed to be $\mathcal{O}(T^{-\frac{p-1}{p}})$.
% For infinite variance, \cite{pmlr-v235-gorbunov24a,pmlr-v202-sadiev23a,nguyen2023high} have shown that, for known iteration number, there exists a constant step size enabling clipped SGD to achieve a convergence rate of $\mathcal{O}(T^{-\frac{p-1}{p}})$. 
However, in some optimization problems, such as object tracking \cite{Shahrampour2016DistributedOO}, iteration number could be unknown. In \cite{nguyen2023improved}, the convergence rate $\mathcal{O}(T^{-\frac{p-1}{p}})$ is achieved by using diminishing step sizes. However, in distributed optimization, the technique in \cite{nguyen2023improved} is not applicable, since the consensus errors among the agents no longer converge.}
  \begin{table*}[!t]
\caption{Comparison among related studies on optimization with noisy gradient.}
\label{table_comparison}
\centering
\resizebox{\linewidth}{!}{
\begin{tabular}{cccccc}
\toprule
Ref. & Loss function &  Problem type & Noise variance & Convergence rate & Requiring prior knowledge(type)\\
\midrule
Ghadimi et al. \cite{doi:10.1137/110848864} & Strongly Convex & Centralized & Finite & $\mathcal{O}(T^{-\frac{1}{2}})$ & No \\
Rogozin et al. \cite{9683110} & Strongly Convex & Distributed & Finite & $\mathcal{O}(T^{-\frac{1}{2}})$ & No \\
Nguyen et al. \cite{nguyen2023high} & Convex & Centralized & Infinite & $\mathcal{O}(T^{-\frac{p-1}{p}})$ & Yes (T) \\
Sadiev et al. \cite{pmlr-v202-sadiev23a} & Convex & Centralized & Infinite & $\mathcal{O}(T^{-\frac{p-1}{p}})$ & Yes (T) \\
Gorbunov et al. \cite{pmlr-v235-gorbunov24a} & Convex & Distributed & Infinite & $\mathcal{O}(T^{-\frac{p-1}{p}})$ & Yes (T) \\
Nguyen et al. \cite{nguyen2023improved} & Convex & Centralized & Infinite & $\mathcal{O}(T^{-\frac{p-1}{p}})$ & No \\
\textbf{This work} & Convex & Distributed & Infinite & $\mathcal{O}(T^{-\frac{p-1}{2p}})$ & No \\
\bottomrule
\end{tabular}}
\end{table*}

 In this paper, distributed optimization is studied via a network of agents, where each agent only
has access to a stochastic gradient of its own objective function, and can communicate with its neighbors via a network.  Different from works \cite{10295561,10103556,li2023convergence}, where the stochastic gradient is assumed to satisfy a lighted-tailed noise model, here we consider the case satisfying the heavy-tailed distribution. Compared with the lighted-tailed noise model, the heavy-tailed  noise model is  mathematically and practically more general. The light-tailed noise model usually guarantees
the boundedness of the variance, which leads to the validity of algorithms'
convergence in high probability \cite{10295561,10103556,li2023convergence}. However, in  the heavy-tailed case, the variance is not necessarily bounded. 
% Hence, heavy-tailed noise is more general and the techniques in \cite{10295561,10103556,li2023convergence} are invalid.
To
solve this problem, a distributed clipped stochastic gradient descent algorithm is proposed. In the algorithm, the clipping strategy is employed to bound the heavy-tailed noise. 
 % This technique is already used in centralized optimization \cite{nguyen2023high,nguyen2023improved,gorbunov2023highprobability,armacki2023highprobability}, which motivate us to study the convergence of distributed stochastic algorithms in high probability. 
%  The advantage of employing the clipping operator is that the clipped gradient is bounded. 
 With the usage of the Markov's inequality, the  high probability bound of the heavy-tailed noise is achieved.
%  In the distributed setting, the difficulty lies behind the coupling of the variables since the network error among agents needs consideration. 
%  Moreover, the local minimizer is no longer the global minimizer in distributed case, which needs more analysis of the bound between the state generated by the iteration and the global minimizer. 
Combining the probability theory, convex
optimization and consensus theory, we prove that if the graph sequence is $B$-strongly connected, then the convergence is ensured with high probability. Moreover, the fastest high probability convergence rate is derived.
 
 This paper is organized as follows. In Section \ref{Sec2}, mathematical preliminaries on heavy-tailed noise and graph theory are introduced. In Section \ref{Sec3}, we
 state our main result and give its proof. In Section \ref{Sec4}, simulation
 examples are presented. Section \ref{Sec5} concludes the whole paper.
 
 \textbf{Notations}: 
 $[n]$ represents the set $\{1,2,\cdots,n\}$ for any integer $n$. $\langle\cdot, \cdot\rangle$ is the standard inner product operator. $\mathbb{E}[\cdot]$ is the expectation operator. $\|\cdot\|$ represents the $2$-norm operator. Let $\mathbb{R}^d$ be the $d$-dimensional real vector space. For differentiable function $f:\mathbb{R}^d\to  \mathbb{R}$, we denote the gradient of $f(x)$ with respect to $x$ by $\nabla f$. For a matrix $W$, $[W]_{ij}$ denotes the matrix entry in the $i^{th}$ row and $j^{th}$ column. 
\section{Preliminaries and Problem Formation}\label{Sec2}
\subsection{Distributed Optimization with Stochastic Gradient Information}
% \textcolor{black}{We focus on the decentralized distributed case.} 
Consider a multi-agent system consisting of $N$ agents, labeled by set  $\mathcal{V}=[N]$. Agents communicate with each other via a time-varying graph sequence  $\mathcal{G}_t$. The goal of agents is to cooperatively solve the following optimization problem:
	\begin{equation}\label{offlineopt}
		\begin{split}
			\min f(x), \text{ with } f(x) =\textcolor{black}{\frac{1}{N}}\sum_{i=1}^{N} f_{i}(x), \text { subject to } x \in \mathbb{R}^d.
		\end{split}
	\end{equation}
 where  $f_{i}: \mathbb{R}^d \rightarrow \mathbb{R}$ \textcolor{black}{is the local objective function.} It is convex and is only available to agent $i$. 
%  $T \in \mathbb{N}$ is unknown to the agents
%  , and  $\mathcal{X} = \mathbb{R}^{d} $. 
%  At each iteration time  $t \in[T]$, agent $i$ selects a state $ x_{i,t} \in \mathcal{X}$. After the state is selected, the information associated with the local cost function  $f_{i}$  is received by agent  $i$, that is, information on cost functions is not available before decisions are made by agents. 
 % In this scenario, at each iteration time  $t$, 

	Here we assume that agent $i$ only has access to a noisy gradient of $f_{i}$ for any $i\in V$, denoted by $\widehat{\nabla}f_{i, t}: \mathbb{R}^d\to  \mathbb{R}^d$. For the estimation of the gradient, we make the following assumptions:
	\begin{assu}\label{assu1}
		For any $x\in \mathcal{X}$,
		\begin{equation*}
			\mathbb{E}[\widehat{\nabla}f_{i}(x)|x]=\nabla f_{i}(x),~\forall~i\in [n].
		\end{equation*}
	\end{assu}
	% \begin{defn}
	% 	A random variable $X$ is $\sigma$-sub-Gaussian if
	% 	\begin{equation*}
	% 		\mathbb{E}[\exp (\lambda^2X^2)]\leq \exp (\lambda^2\sigma^2)
	% 	\end{equation*}
	% 	for all $\lambda$ such that $|\lambda|\leq \frac{1}{\sigma}$.
	% \end{defn}
	
	In our work, we propose a high probability convergence in the heavy-tailed noise regime. More specifically, we analyze the heavy-tailed noise
 model proposed by \textcolor{black}{\cite{Yudin1983ProblemCA} and \cite{zhang2020adaptive}}, where the gradient noise has bounded $p$-th moment.
	\begin{assu}\label{assu2}
		For any $x\in \mathbb{R}^d$, $i\in [N]$, 
		\begin{align*}
		   \mathbb{E}[\|\widehat{\nabla}f_{i}(x)-\nabla f_{i}(x) \|^p|x]\le \sigma^p 
		\end{align*}
		for some $\sigma>0$ and $p\in (1,2]$.
	\end{assu}
	
	This is commonly referred to as heavy-tailed noise, as opposed to light-tailed noise such as those that are distributed according to sub-Gaussian distributions. Mathematically, heavy-tailed noise is more general in stochastic optimization. The  variance of the light-tailed noise, e.g. sub-Gaussian noise, is finite since its tail decays faster than Gaussian distribution. Thus, Assumption \ref{assu2} can be derived in the light-tailed settings but the opposite is not true.
	% Assumptions \ref{assu1} and \ref{assu2} are also used in \cite{lu2023highpro, liu2023high}, where the convergence of the offline algorithms in high probability is studied. Different from them, we will study the online case. 
		% \begin{assu}\label{assu3}
		% $\mathcal{X}$ is convex and compact.
		% 	\end{assu}
	\begin{assu}\label{assu3}
		For any $i\in [N]$, $f_{i}$ is $L$-smooth, i.e., its gradient is $L$-Lipschitz. For all $x,y\in \mathbb{R}^d$,
		\begin{equation*}
			|\nabla f_{i}(x)- \nabla f_{i}(y)|\leq L \| x-y\|.
		\end{equation*}
   Moreover, if $f_{i}$ is $L$-smooth, it has a quadratic upper bound, i.e., For all $x,y\in \mathbb{R}^d$,
        \begin{equation*}
		 f_{i}(y)\le f_{i}(x) +\langle \nabla f_{i}(x), y-x \rangle +\frac{L}{2}\|y-x\|^2.
		\end{equation*}
	\end{assu}

	Agents communicate with each other via an directed graph $\mathcal{G}_t=(\mathcal{V}, W_t, \mathcal{E}_t)$, where $\mathcal{V}$ denotes the set of agents, $\mathcal{E}_t$ denotes the set of edges where the elements are denoted as $(i,j)$ if and only if agent $i$ can
receive a message from agent $j$ at time $t$, and $W_t=([W_t]_{ij})_{n\times n}$ denotes the weighted matrix. We denote the set of
incoming neighbors of agent $i$ at time $t$ by
\begin{equation*}
    \mathcal{N}_{i,t}=\{j|(i,j)\in \mathcal{E}_t \}\cup \{i\}.
\end{equation*}
For a fixed topology $\mathcal{G}=(\mathcal{V}, W, \mathcal{E})$, a path of length $r$ from $i_1$ to $i_{r+1}$ is a sequence of $r+1$ distinct nodes $i_1, \cdots, i_{r+1}$ such that $(i_{q},i_{q+1})\in \mathcal{E}$, for $q\in [r]$. If there exists a path between any two nodes, then $\{\mathcal{G}\}$ is said to be strongly connected. For $\{\mathcal{G}_t\}$, an $B$-edge set is defined as $\mathcal{E}_{B,t}=\bigcup_{l=0,\cdots,B-1}\mathcal{E}_{t+l}$ for some constant $B>0$. We call that $\{\mathcal{G}_t\}$ is $B$-strongly connected if the directed graph with vertex $\mathcal{V}$ and edge $\mathcal{E}_{B,t}$ is strongly connected for any $t\ge 1$.
	\begin{assu}\label{assu6}
		For all $t \ge 1$, the weighted graphs $\mathcal{G}_t=(\mathcal{V}, W_t, \mathcal{E}_t)$ satisfy:\\
  (a) There exists a scalar $\eta \in (0, 1)$ such that $[W_t]_{ij}\ge \eta  $ if $j\in \mathcal{N}_{i,t}$. Otherwise, $[W_t]_{ij}=0 $.\\
  (b) The weighted matrix is doubly stochastic, i.e.,
		\begin{equation*}
			\sum_{i=1}^{n}[W_t]_{ij}=\sum_{j=1}^{n}[W_t]_{ij}=1.
		\end{equation*}\\
  (c) $\{\mathcal{G}_t\}$ is $B$-strongly connected. That is, there exists a scalar $B > 0$ such that the graph
$(\mathcal{V},\mathcal{E}_{B,t})$ is strongly connected for any $t \ge 1$.
	\end{assu}
  For any $t \ge s$, we denote
\begin{align*}
   \left\{\begin{array}{l}
   \Phi(t, s)=W_{t-1} \cdots W_{s+1} W_{s}, \text { if } t>s \\
\Phi(t, s)=I_{n}, \text { if } t=s.
\end{array}\right. 
\end{align*}
\textcolor{black}{
Under the assumption above, we have the following lemma \cite{sundhar2012new}.
\begin{lemma}[\cite{sundhar2012new}]\label{lem_matrix}
    Under Assumption \ref{assu6}, for any $k\ge s\ge 1$, 
\begin{align*}
    \left|[\Phi(k, s)]_{i, j}-\frac{1}{N}\right| \leq \gamma \beta^{k-s},
\end{align*}
where
\begin{align}\label{gamma_eta}
    \gamma=\left(1-\frac{\eta}{4 N^{2}}\right)^{-2} \quad \beta=\left(1-\frac{\eta}{4 N^{2}}\right)^{\frac{1}{B}}.
\end{align}
\end{lemma}}
\subsection{Distributed Algorithm with Clipping Operator}
To solve problem (\ref{offlineopt}), we provide our algorithm, which is called the distributed clipped stochastic gradient descent
\begin{equation}\label{algorithmoffline}
\begin{cases}
    y_{i,t}=\sum_{j=1}^{N}\left[W_{t}\right]_{i j}x_{j,t}\\
    \widetilde{\nabla} f_{i}(x_{i,t})=\text{clip}_{\lambda_t}(\widehat{\nabla} f_{i}(x_{i,t}))\\
    x_{i,t+1}=y_{i,t}-\eta_t\widetilde{\nabla} f_{i}(x_{i,t})
\end{cases}
\end{equation}
where $\{\eta_{t}\}_{t=1}^{T}$ is the step size sequence, $x_{i,t}\in \mathbb{R}^d$ represents the state of agent $i$ with initial state $x_{i,1}\in \mathbb{R}^d$. $\widehat{\nabla}f_{i}: \mathbb{R}^d\to \mathbb{R}^d$ is the non-biased estimation of the gradient ${\nabla}f_{i}$.
	The step size sequence is non increasing and positive and the clipping operator $\text{clip}_{\lambda_t}(\cdot): \mathbb{R}^d\to \mathbb{R}^d$ is 
 \begin{align}\label{clipping_operator}
     \text{clip}_{\lambda_t}(y)=\min \{1,\frac{\lambda_t}{\|y\|}\}y
 \end{align}
 \textcolor{black}{where $\{\lambda_{t}\}_{t=1}^{T}$ is the clipping parameter sequence.}
 \begin{rem}
  Algorithm (\ref{algorithmoffline}) is designed by combining the consensus algorithm \cite{4631522,sundhar2012new,5400096} and the clipping operator\textcolor{black}{\cite{Pascanu2013,nguyen2023high,nguyen2023improved,gorbunov2020stochastic,zhang2020adaptive}}. But works \cite{nguyen2023high,nguyen2023improved,gorbunov2020stochastic,zhang2020adaptive} are conducted by centralized approaches, and the results are not applicable to distributed cases due to the coupling of the variables. Algorithm (\ref{algorithmoffline}) runs by using the clipped unbiased gradient estimator of the local cost function and the state information received from its neighbors. Therefore, algorithm (\ref{algorithmoffline}) is distributed.
 \end{rem}

\section{Main Results}\label{Sec3}
In this section, we will state our main results. First we analyze the network error, i.e., the error between each agent’s state
 and their average value at each iteration under (\ref{algorithmoffline}). 
Note that $x_{i,t}$ can be represented by the iteration in Lemma \ref{lem_matrix}. By (\ref{gamma_eta}) and algorithm (\ref{algorithmoffline}), we can bound the network error as follows
\begin{lemma}\label{lem_networkerr}
    Under Assumption \ref{assu6}, by algorithm (\ref{algorithmoffline}), for all $ i \in [N] $ and $ t \geq 1$ \textcolor{black}
{\begin{align*}
	\left\|x_{i, t+1}-\overline{x}_{t+1}\right\|& \leq  N\gamma \beta^{t} R_1+
 N\gamma \sum_{l=1}^{t}\beta^{t-l}\lambda_{l}\eta_l,
\end{align*}}
% and
% \begin{align*}
% 	\sum_{i=1}^{N}\left\|\lambda_{i, t+1}-\overline{\lambda_t}_{t+1}\right\|& \leq  N \gamma \beta^{t} B_\lambda_t+2 NB_y\eta_{t}\\
%  &+\gamma N^2 B_y \sum_{l=1}^{t-1}\beta^{t-l}\eta_l,
% \end{align*}
where $\bar{x}_t=\frac{1}{N}\sum_{i=1}^{N}x_{i,t} $, $ \gamma $ and $ \beta $ are defined in (\ref{gamma_eta}) and $R_1=\max_{i}\|x_{i,1}\|$.
\end{lemma}

\textcolor{black}{\textit{Proof:} See Appendix \ref{proof_networkerr}.}

$\gamma$ and $\beta$ are significant factors that influence the error bound. If the connected period $B$ increases, the error bound decreases. After the analysis of the network error, the high probability convergence of the distributed optimization can be bounded. 
% In the following part we provide main results.
% \subsection{High Probability Convergence}
% In this part, the high probability convergence is derived. 
The main difficulty arises when considering the deviation between the exact gradient and the clipped one, which is not unbiased. To bound the error induced by the heavy-tailed noise, we propose a method that can bound the state $\|x_{i,t}-x^*\|$ by $\|x_{i,1}-x^*\|$, where the coupling of the variables adds the difficulty since $x^*\ne \operatorname{argmin} f_{i}(x)$. First, we analyze the bound of $\|\bar{x}_t-x^*\|$. Here we \textcolor{black}{denote $\theta_{i,t}=\widetilde\nabla f_{i}(x_{i,t})-\nabla f_{i}(x_{i,t})$} and highlight the appearance of the term $\langle \theta_{i,t},\bar{x}_{t}-x^{\star}\rangle=\langle\widetilde\nabla f_{i}(x_{i,t})-\nabla f_{i}(x_{i,t}),\bar{x}_{t}-x^{\star}    \rangle$.
\begin{lemma}\label{lem_boundofaverage}
    Under Assumptions \ref{assu1}-\ref{assu3}, by algorithm (\ref{algorithmoffline}),
    \textcolor{black}{\begin{equation}\label{boundofaverage}
        \begin{aligned}
            &\eta_t (f(\bar{x}_t)-f(x^*)) \le \frac{\eta_t^2\lambda_t^2}{2}\\
    &+\frac{1}{2}(\|\bar{x}_{t}-x^{*}\|^2-\|\bar{x}_{t+1}-x^{*}\|^2)\\
    & + \frac{1}{N}\sum_{i=1}^{N} \eta_t \langle\theta_{i, t}, \bar{x}_{t}-x^{*}\rangle +\frac{L\eta_t}{2} \| \bar{x}_{t}-x_{i,t} \|^{2}.
        \end{aligned}
    \end{equation}}
\end{lemma}

\textit{Proof: } See Appendix \ref{proofboundofaverage}.

Note that $x^*=\operatorname{argmin} f(x)$, the left hand side of (\ref{boundofaverage}) is non-negative. We have access to derive the bound of $\|\bar{x}_{n+1}-x^*\|^{2}$ by $\|\bar{x}_{i}-x^*\|$, $i\in [n]$ inductively. With the help of Lemma \ref{lem_networkerr}, we can bound $\|\bar{x}_{n+1}-x^*\|^{2}$ via the step size $\eta_t$ and $\lambda_t$. Here we denote $z_t=\frac{1}{a_0+4\max_{i\le t}\|\bar{x}_{i}-x^*\|}$ with $a_0=\max \{1,\log \frac{1}{\delta}\}$.
 \begin{lemma}\label{lem_statebound}
     Under Assumptions \ref{assu1}-\ref{assu3}, by algorithm (\ref{algorithmoffline}), for any $n\in [T]$, 
\begin{align*}
    &z_n\|\bar{x}_{n+1}-x^*\|^2\le z_1\|\bar{x}_{1}-x^*\|^2 + \textcolor{black}{\frac{4L\eta_1z_1}{1-\beta}N^2\gamma^2R_1^2} \\
    &+ \textcolor{black}{4N^2\gamma^2(\frac{1}{1-\beta})^2}Lz_1\sum_{t=1}^{n}\eta_t^3\lambda_t^2+z_1\sum_{t=1}^{n}\eta_t^2\lambda_t^2\\
    &+ \sum_{t=1}^{n}\frac{1}{N} \sum_{i=1}^{N} 2z_t\eta_t \langle\theta_{i, t}, \bar{x}_{t}-x^{*}\rangle.
\end{align*}
 \end{lemma}
 
 \textit{Proof: } See Appendix \ref{proofstate}.

Before obtaining the bound, we need to handle the term $z_t\eta_t\langle \theta_{i,t},\bar{x}_{t}-x^{\star}\rangle$, the difficulty is that the noise is heavy-tailed. To handle the difficulty, we refer to the lemma presented in \cite{nguyen2023high}. Let
\begin{align*}
    &\theta_{i,t}^u=\widetilde{\nabla}f_{i}(x_{i,t})-\mathbb{E}[\widetilde{\nabla}f_{i}(x_{i,t})|\mathcal{F}_t],\\&
    \theta_{i,t}^b=\mathbb{E}[\widetilde{\nabla}f_{i}(x_{i,t})|\mathcal{F}_t]-\nabla f_{i}(x_{i,t})
\end{align*}
where $\mathcal{F}_t=\sigma(\widehat{\nabla} f_{1},\cdots,\widehat{\nabla} f_{t-1})$ denotes the $\sigma$-field generated by the unbiased estimator by $t$ ($(\mathcal{F}_t)_{t\ge 0}$ is also known as the natural filtration).
Note that $\theta_{i,t}=\theta_{i,t}^u+\theta_{i,t}^b$ and $\theta_{i,t}^u$ is unbiased under $\mathcal{F}_t$. With the notation above, we provide the following lemma. 
\begin{lemma}[\textcolor{black}{\cite{pmlr-v202-sadiev23a}}]\label{lem_heavytaied}
    Under Assumptions \ref{assu1} and \ref{assu2}, for $t\ge1$ and $i\in [N]$,  we have
    \begin{align*}
        \|\theta_{i,t}^u\|\le 2\lambda_t.
    \end{align*}
    Furthermore, if $\|\nabla f_{i}(x_{i,t})\|\le \frac{\lambda_t}{2}$, then
    \begin{align*}
        &\|\theta_{i,t}^b\|\le 4\sigma^p\lambda_{t}^{1-p},\\
        &\mathbb{E}[\|\theta_{i,t}^u\|^2|\mathcal{F}_t]\le 16\sigma^p\lambda_t^{2-p}.
    \end{align*}
\end{lemma}

By the $L$-smoothness of the gradient of the $f_{i}$, $i\in [N]$, we have access to bound the exact gradient $\|\nabla f_{i}(x_{i,t})\|$ by
\begin{equation}\label{bound_extgradient}
    \begin{aligned}
     &\|\nabla f_{i}(x_{i,t})\|\\
    &\le\|\nabla f_{i}(x_{i,t})-\nabla f_{i}(x^*)\|+\|\nabla f_{i}(x^*)\|\\
    &\le L\|x_{i,t}-x^*\|+\|\nabla f_{i}(x^*)\|\\
    &\le L\|x_{i,t}-\bar{x}_t\|+L\|\bar{x}_t-x^*\|+\|\nabla f_{i}(x^*)\|
    \end{aligned}
\end{equation}
which implies that $\|\nabla f_{i}(x_{i,t})\|$ can be bounded via appropriate choice of $\lambda_t$ once the high probability bound of $\|\bar{x}_t-x^*\|$ is achieved. The outline of the proof is as follows. First, by Lemma \ref{lem_statebound}, and with the help of the decreasing coefficient $z_t$, $z_t\eta_t\langle \theta_{i,t}^u,\bar{x}_{t}-x^{\star}\rangle$ can be bounded by its quadratic conditional expectation without considering $\|\bar{x}_n-x^*\|$ for any $n\in [t-1]$. Second, by the induction on the number of iterations and (\ref{bound_extgradient}) we have access to bound the error induced by the noise by Lemma \ref{lem_heavytaied}. Last, combining all bounds together, $\|\bar{x}_t-x^*\|$ can be bounded by  $\|\bar{x}_1-x^*\|$.

% Summing inequality (\ref{main_error}) by $[n]$, we can bound $\|x_{i, n+1}-x^{\star}\|^2$ since the left hand of inequality (\ref{main_error}) is positive, i.e.,
% \begin{align*}
%     &\sum_{t=1}^{n}\eta_t (f_{i, t}\left(x_{i, t}\right)-f_{i, t}\left(x^{\star}\right ))\\
%     &\ge \sum_{t=1}^{n}\eta_n (f_{i, t}\left(x_{i, t}\right)-f_{i, t}\left(x^{\star}\right )) \ge 0.
% \end{align*}
% Thus, it yields that
% \begin{align*}
%   &\|x_{i, n+1}-x^{\star}\|^2 \le \|x_{i, 1}-x^{\star}\|^2 + \sum_{t=1}^{n} (11+4N\gamma \frac{1}{1-\beta}) \lambda_t^2\eta_t^2\\& + 4NR_1\gamma \lambda_1\eta_1\frac{1}{1-\beta}+ 2\sum_{t=1}^{n} \eta_t\langle \theta_{i,t},x_{i, t}-x^{\star}\rangle.
% \end{align*}
% The inequality above is satisfied by the Fubini's theorem
% \begin{align*}
%     &\sum_{t=1}^{T}\lambda_t\eta_t\sum_{l=1}^{t-1}\beta^{t-l}\lambda_l\eta_l = \sum_{l=1}^{T-1}\sum_{t=l+1}^{T-1}\beta^{t-l}\lambda_t\eta_t\lambda_l\eta_l\\
%     &=\sum_{l=1}^{T-1}\frac{1}{1-\beta} \lambda_{l+1}\eta_{l+1}\lambda_l\eta_l
%     \le \sum_{l=1}^{T}\frac{1}{1-\beta} \lambda_l^2\eta_l^2.
% \end{align*}
Now we provide the appropriate choice of $\eta_t$ and $\lambda_t$. Denote $\Delta_{n}=\|\bar{x}_{n}-x^{\star}\|$, \textcolor{black}{$A=4LN^2\gamma^2(\frac{1}{1-\beta})^2$ and $E=\frac{4L}{1-\beta}N^2\gamma^2R_1^2$}.
 For fixed $\delta>0$, we denote $a_0=\max \{1,\log \frac{1}{\delta}\}$ and adapt the clipping parameter $\lambda_t=\lambda t^\alpha$, $\alpha> 0$ and step sizes $\eta_t=\frac{1}{m_{t}t^{\kappa}}$, \textcolor{black}{where  $m_{t}=m(\log t + b_1)^{2}$, $m\ge 1$ and $b_1\ge 1$.} $\lambda_t$ and $\eta_t$ satisfy the following condition
\begin{align}\label{condition}
  \begin{cases}
    % AC_{\kappa}\lambda_t^2+B\lambda_t \le \frac{d^2\Delta_{i,1}}{2}\\
    % \lambda_t \ge 4N\sqrt{L\Delta_{i,1}}
    \kappa\ge \max\{\alpha+\frac{1}{2},1-(p-1)\alpha \}\\
    \textcolor{black}{AC_{2\alpha,3\kappa}\lambda^2+(C_{2\alpha,2\kappa}N+Em)m\le m^3\Delta_{1}^2}\\
    \lambda \ge 2L(9\Delta_1+5a_0)+2LD+2B^*\\
    \textcolor{black}{m}\ge \Delta_{1}^{-1}\lambda^{1-p}\sigma^p C_{(1-p)\alpha,\kappa} \\
    \textcolor{black}{m}\ge 6\Delta_1^{-0.5}\lambda^{1-p/2}\sigma^{p/2}C_{2\alpha,2\kappa}^{0.5}
\end{cases}  
\end{align}
where \textcolor{black}{$C_{a,b}=\sum_{t=1}^{\infty}\frac{1}{t^{b-a}(\log t + 1)^{2}}$} and $B^*$ is the largest norm of the gradient $\|\nabla f_{i}(x^*)\|$, $i\in [N]$. The reason of choosing $\lambda_t$ that increases to infinity is to bound $\sum_{t=1}^{n} z_t\eta_t\langle \theta_{i,t}^b,\bar{x}_t-x^{*}\rangle$ by Lemma \ref{lem_heavytaied}. It is worth highlighting that the network error $\|x_{i,t}-\bar{x}_t\|$ is bounded by \textcolor{black}{$D=N\gamma R_1+
 N\gamma\frac{\lambda}{(1-\beta)m}$}.
\begin{rem}
    Parameters $\alpha$ and $\kappa$ satisfying condition (\ref{condition}) are easy to be selected. For instance, \textcolor{black}{$(\alpha,\kappa)=(\frac{1}{3p},\frac{2p+1}{3p})$}. Moreover, the second inequality of the condition (\ref{condition}) involves the global information such as the connectivity period. 
    In practical application, we can select \textcolor{black}{$m$} as large as possible without knowing the global information.
\end{rem}

Then we provide the high probability bound of the term $\|\bar{x}_t-x^*\|$. 
% \textcolor{black}{Motivated by \cite{nguyen2023improved}, we have the following lemma. }
\begin{lemma}\label{lem_stateboundbynoise}
Under Assumptions \ref{assu1}-\ref{assu6}, for $\eta_t$ and $\lambda_t$ satisfying condition (\ref{condition}), then, for any $n\ge 1$ , with probability at least $1-\delta$, 
\begin{align*}
    \sum_{t=1}^{n} 2z_t\eta_t\langle \theta_{i,t},\bar{x}_{ t}-x^{\star}\rangle \le \frac{5\Delta_{1}}{4}+\log \frac{1}{\delta}.
\end{align*}
% Under Assumptions \ref{assu1}-\ref{assu6}, for $\eta_t$ and $\lambda_t$ satisfying condition (\ref{condition}), $E_n$ happens with probability at least $1-\delta$.
\end{lemma}

\textit{Proof:} See Appendix \ref{proofnoise}.

% Now the $\Delta_{t}$ is bounded by induction, thus Lemma \ref{lem_heavytaied} is feasible. We gain the high probability convergence as follows
Now we present our main result.
\begin{thm}\label{thm_highprobconvergence}
  Under Assumptions \ref{assu1}-\ref{assu6},  if $\eta_t$ and $\lambda_t$ satisfy condition (\ref{condition}), then by algorithm (\ref{algorithmoffline}), with probability at least $1-\delta$,
  \textcolor{black}{\begin{align*}
     \sum_{t=1}^{T} (f(\bar{x}_t)-f(x^*)) \le \frac{9}{2\eta_T}\Delta_{1}^2+ \frac{1}{\eta_T} \textcolor{black}{\log^2 \frac{1}{\delta}}.
  \end{align*}}
\end{thm}

\textcolor{black}{\textit{Proof:} See Appendix \ref{proof_thm}.}

Moreover, the fastest high probability convergence rate induced by the algorithm (\ref{algorithmoffline}), is
\begin{corollary}
    Under the the same conditions
stated in Theorem \ref{thm_highprobconvergence}, if $\lambda_t=\lambda t^{\frac{1}{2p}}$ and $\eta_t=\frac{1}{m_{t}t^{\frac{1}{2}+\frac{1}{2p}}}$, then with probability at least $1-\delta$,
  \begin{align*}
   \frac{ \sum_{t=1}^{T} f(\bar{x}_t)-f(x^*)}{T} \textcolor{black}{\le} \mathcal{O}(\textcolor{black}{(N \rho^{\frac{5}{3}}+B^{*}\rho^{\frac{2}{3}})}T^{\frac{1}{2p}-\frac{1}{2}}\log^2 (\frac{T}{\delta})).
  \end{align*}
\textcolor{black}{where $B^*=\max_{i\in [N]}\|\nabla f_{i}(x^*)\|$ and $\rho=\frac{1}{1-\beta}$.}
\end{corollary}

\textit{Proof:} By the first line in condition (\ref{condition}), $\kappa\ge \max\{\alpha+\frac{1}{2},1-(p-1)\alpha \}$. Note that $\min\max\{\alpha+\frac{1}{2},1-(p-1)\alpha \}=\frac{1}{2p}+\frac{1}{2}$ where $\alpha=\frac{1}{2p}$. \textcolor{black}{By the second and third lines in (\ref{condition}) we know that $m\sim\frac{1}{(1-\beta)^{2/3}}\lambda^{2/3}$ and $\lambda\sim\mathcal{O}(\frac{N}{1-\beta}+B^{*})^{1.5}$.} Substituting the minimum into Theorem \ref{thm_highprobconvergence} we obtain the fastest high probability convergence rate.

\begin{rem}
% \textcolor{black}{Since $p\in (1,2]$, the convergence rate indicates a linear speedup with respected to $\sigma$.} 
If the deviation term $\widehat{\nabla}f_{i}(x)-\nabla f_{i}(x) $ has a finite variance, algorithm (\ref{algorithmoffline}) ensures the fastest convergence rate $\mathcal{O}(T^{-\frac{1}{4}} \log^2 (\frac{T}{\delta}))$. 
\end{rem}

\section{Simulation}\label{Sec4}
\begin{figure}
     \centering
     \includegraphics[scale=0.25]{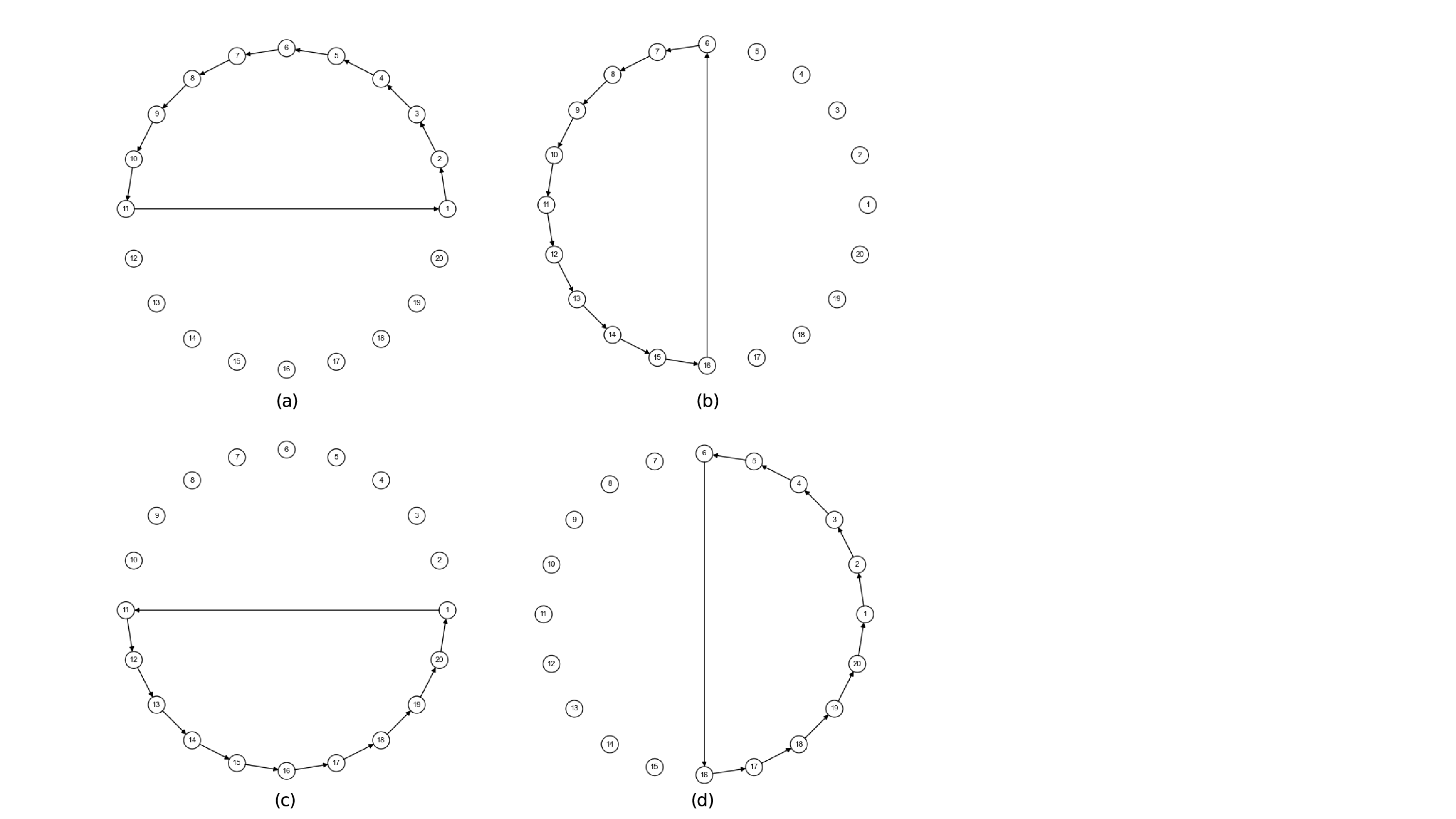}
     \caption{The Time-varying Graph}
     \label{fig:time-varying}
 \end{figure}

\textcolor{black}{In this section, we use our methods to solve a distributed classification problem. We consider $N=20$ agents embedded in a directed and time-varying graph. Each sensor communicates with its neighbors via a time-varying graph shown in Fig.\ref{fig:time-varying}. }
\textcolor{black}{The weight
 of each edge in Fig.\ref{fig:time-varying} is assumed to be $0.8$. The switching order is given by $(a)\to (b)\to (c)\to (d)\to (a)\to \cdots$. Note that
 the union of the possible graphs is strongly connected. Then,
 the connectivity of the graph in Fig.\ref{fig:time-varying} satisfies conditions in
 Assumption \ref{assu6} with $B = 4$.}

\textcolor{black}{Each agent $i$ possesses a local dataset $\mathcal{D}_i$, comprising $m_i$ training samples $\eta_{ij}=\{x_{ij}^{\top},y_{ij}\}_{j=1}^{m_{i}}$. The $j$-th sample at node $i$ is a tuple $\{x_{ij}^{\top},y_{ij}\}\subset \mathbb{R}^{d}\times \mathbb{Z} $, where $x_{ij}^{\top}$ represents the }\textcolor{black}{features and $y_{ij}\in \mathbb{Z}$ represents the
 label. The total number of labeled data points is $m=\sum_{i=1}^{N}m_{i}$.
We choose the average of the local logistic regression
 loss function:
\begin{equation}\label{simuprob}
    \begin{aligned}
        &f(\omega)=\frac{1}{N} \sum^{N}_{i=1}f_{i}(\omega)\\
        &=\frac{1}{N} \sum^{N}_{i=1} \frac{1}{m_{i}} \sum^{m_{i}}_{j=1} \frac{1}{2} \|x_{ij}\omega-y_{ij}\|^2+\frac{\lambda}{2}\|\omega\|^2.
    \end{aligned}
\end{equation}
We use 2000 samples for training with the regularization parameter $\lambda=0.1$. Note that the local gradient is not accurate, i.e.,
	\begin{equation*}
		\widehat{\nabla}f_{i, t}=\frac{1}{m_{i}} \sum^{m_{i}}_{j=1} x_{ij}^{\top}x_{ij}\omega-y_{ij}+\lambda\omega+ 0.2\xi_{i,t}
	\end{equation*}
	where $\xi_{i,t}\in \mathbb{R}^{d}$ is the heavy-tailed noise. We consider a typical noise, the t-distribution $t_2$ with the probability density function
	\begin{align*}
	    t_2(x)=\frac{\Gamma(\frac{3}{2})}{\Gamma(\frac{1}{2})
	   \sqrt{2\pi}}(1+\frac{x^2}{2})^{-\frac{3}{2}}
	\end{align*}
	where $\Gamma(x)=\int_{0}^{
	+\infty}t^{x-1}e^{-t}dt$. We can verify that $t_2$ has zero expectation and unbounded variance. 
	Then we solve problem (\ref{simuprob}) by our algorithm. In simulation, we choose $T=2000$. The initial values are randomly chosen from $\mathbb{R}^{d}$. We perform the algorithm on real-world datasets and compare it with the distributed stochastic gradient descent.}

\begin{figure}
    \centering
    \includegraphics[scale=0.5]{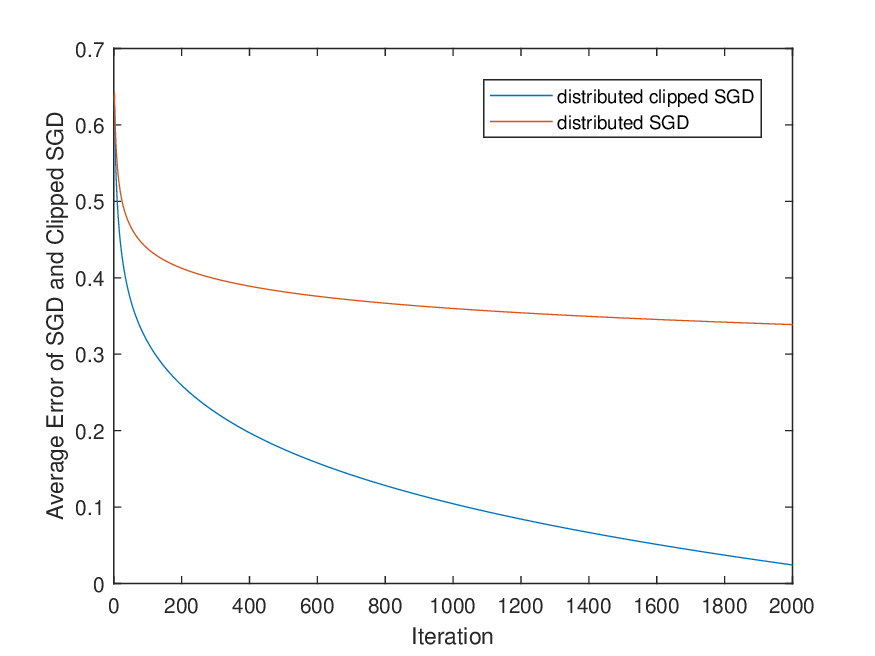}
    \caption{The Average Convergence of the Distributed Clipped SGD and the Distributed SGD on a9a}
    \label{fig:a9a}
\end{figure}

\textbf{Example 1.} \textcolor{black}{We perform classification on the a9a dataset, where the label $y_{ij}\in \{0,1\}$ , $d=1558$. We choose the parameter
\begin{align*}
	    \eta_t&=\frac{5}{t^{0.75}(1 +\log t)^{2}},\\
	    \lambda_t&=2t^{0.25}.
\end{align*}
We run the algorithm in a single round and compare the performance to the traditional distributed SGD algorithm in }\textcolor{black}{\cite{10295561} for our problem with the same step size. The average convergence $\frac{1}{T}\sum_{t=1}^{T}f(\bar{\omega}_{t})-f(\omega^{*}) $ under the algorithms is depicted in Fig.\ref{fig:a9a}.}

\textcolor{black}{Note that the algorithm (\ref{algorithmoffline}) admits a lower average error since the error induced by the noise is clipped by the operator.}

\begin{figure}
    \centering
    \includegraphics[scale=0.5]{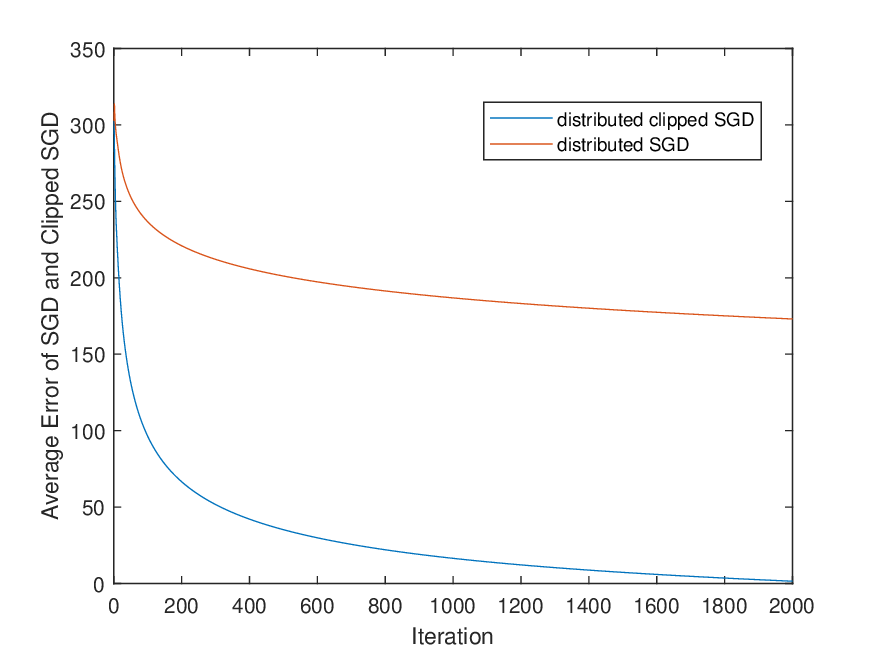}
    \caption{The Average Convergence of the Distributed Clipped SGD and the Distributed SGD on MNIST}
    \label{fig:MNIST}
\end{figure}

\textbf{Example 2.} \textcolor{black}{Now we perform the classification on the MNIST dataset, where the label $y_{ij}\in [10]$ and $d=784$. The parameters are given by
\begin{align*}
	    \eta_t&=\frac{1}{35t^{0.75}(4 +\log t)^{2}},\\
	    \lambda_t&=100t^{0.25}.
\end{align*}
Similarly, we run the algorithm in a single round and compare the performance to the traditional distributed SGD algorithm in \cite{10295561} for our problem with the same step size. The average convergence of the algorithms mentioned is depicted in Fig.\ref{fig:MNIST}.}

\textcolor{black}{Note that the algorithm (\ref{algorithmoffline}) admits a lower average error since the error induced by the noise is clipped by the operator. These observations are consistent with our obtained results.}
\section{Conclusion}\label{Sec5}
%  In this paper, the problem of distributed optimization has been studied by employing
%  a network of agents, where each agent can only communicate
%  with its immediate neighbors via a time-varying directed
%  graph.
 In this paper, we have proposed a distributed clipped stochastic gradient descent algorithm for
 the distributed optimization problem.  
 By implementing the algorithm,
 every agent adjusts its state value by 
 the clipped gradient estimation from its local cost function and the local
 state information received from its immediate neighbors.
 We show that
if the time-varying graph sequence is $B$-strongly connected, the high probability convergence rate is $\mathcal{O}(T^{\frac{1}{2p}-\frac{1}{2}}\log^2 (\frac{T}{\delta}))$. A
 simulation example has been presented to demonstrate the
 effectiveness of our theoretical results. 
 
 When constraints exist in distributed optimization problems,
 more complicated distributed stochastic algorithms, such as the
 one based on the primal-dual strategy, are needed. How to study the convergence of those distributed algorithms in high
 probability is still an interesting open topic. Some other issues may also be considered, such as the case
 with time delays and communication bandwidth constraints,
 which will bring new challenges in the study of convergence in high probability.
\section{Appendix}
\subsection{Proof of Lemma \ref{lem_networkerr}}\label{proof_networkerr}
\textit{Proof of Lemma \ref{lem_networkerr}:} 
Letting $e_{i, t}=x_{i, t+1}-y_{i, t}$, by algorithm (\ref{algorithmoffline}), we have $x_{i, t+1}=\sum_{j=1}^{N}[W_t]_{i j} x_{j, t}+e_{i, t}$. Due to the fact that $W_t$ is doubly stochastic, we obtain 
	\begin{equation}\label{average_x}
		\bar{x}_{t+1}:=\frac{1}{N} \sum_{i=1}^{N} x_{i, t+1}=\bar{x}_{t}+ \bar{e}_{t}=\sum_{\tau=1}^{t} \bar{e}_{\tau}+\bar{x}_1
	\end{equation}
	where  $\bar{e}_{t}=\frac{1}{N} \sum_{i=1}^{N} e_{i, t} $. Meanwhile, letting
	\begin{align*}
		x_{t}:=\left[x_{1, t}^{\top}, x_{2, t}^{\top}, \ldots, x_{N, t}^{\top}\right]^{\top}\\ e_{t}:=\left[e_{1, t}^{\top}, e_{2, t}^{\top}, \ldots, e_{N, t}^{\top}\right]^{\top}
	\end{align*}
	and using  $\otimes$  to denote the Kronecker product, we have $x_{t+1}=(W_t \otimes I_N) x_{t}+\left(I_{N} \otimes I_n\right) e_{t}$, which implies
	\begin{align*}
		&x_{i, t+1}=\sum_{j=1}^{N}\left[\Phi(t,1)\right]_{i j}x_{j,1}\\
  &+\sum_{\tau=1}^{t} \sum_{j=1}^{N}\left[\Phi(t-\tau+1,1)\right]_{i j}  e_{j, \tau}.
	\end{align*}
	Together with (\ref{average_x}), there holds
	\begin{align*}
		x_{i, t+1}-\bar{x}_{t+1}&=\sum_{\tau=1}^{t} \sum_{j=1}^{N}\left(\left[\Phi(t-\tau+1,1)\right]_{i j}-\frac{1}{N}\right)  e_{j, \tau}\\
		+&\sum_{j=1}^{N}\left(\left[\Phi(t+1,1)\right]_{i j}-\frac{1}{N}\right) (x_{j,1}-\bar{x}_1).
	\end{align*}
	It implies
	\begin{equation*}
		\begin{split}
			&\| x_{i, t+1}-\bar{x}_{t+1} \| 
			\leq \sum_{j=1}^{N} \left | \left([\Phi(t+1,1)]_{i j}-\frac{1}{N}\right) \right |\|x_{j,1}-\bar{x}_1\|  \\
			&+\sum_{\tau=1}^{t} \sum_{j=1}^{N} \left | \left([\Phi(t-\tau+1,1)]_{i j}-\frac{1}{N}\right) \right |\eta_{\tau} \lambda_{\tau},
		\end{split}
	\end{equation*}
	where the bound of $\|e_{j, \tau}\|$ is evaluated in \ref{clipping_operator}. By Assumption \ref{assu6}, and using the result in Lemma \ref{lem_matrix}, we have the desired result.
\subsection{Proof of Lemma \ref{lem_boundofaverage}}\label{proofboundofaverage}
\textit{Proof of Lemma \ref{lem_boundofaverage}:} Due to the $L$-smoothness  and the convexity of $f_{i}$, it yields that
\begin{align*}
&f_{i}\left(\bar{x}_{t}\right)-f_{i}\left(x_{i, t}\right) \le\langle\nabla f_{i}\left(x_{i, t}\right), \bar{x}_{t}-x_{i, t}\rangle+\frac{L}{2}\left\|\bar{x}_{t}-x_{i, t}\right\|^{2} \\
&f_{i}\left(x_{i, t}\right)-f_{i}(x^{*}) \le\langle\nabla f_{i}\left(x_{i, t}\right), x_{i,t}-x^{*}\rangle.
\end{align*}
Thus, we have that
\begin{align*}
    &f_{i}\left(\bar{x}_{t}\right)-f_{i}(x^*) \le\langle\nabla f_{i}\left(x_{i, t}\right), \bar{x}_{t}-x^{*}\rangle+\frac{L}{2} \| \bar{x}_{t}-x_{i, t} \|^{2} \\
&=\langle\widetilde{\nabla} f_{i}\left(x_{i, t}\right), \bar{x}_{t}-x^{*}\rangle+\langle\theta_{i, t}, \bar{x}_{t}-x^{*}\rangle+\frac{L}{2} \| \bar{x}_{t}-x_{i,t} \|^{2}
\end{align*}
where $\langle \theta_{i,t},\bar{x}_{t}-x^{\star}\rangle=\langle\widetilde\nabla f_{i}(x_{i,t})-\nabla f_{i}(x_{i,t}),\bar{x}_{t}-x^{\star}    \rangle$. Note that by algorithm (\ref{algorithmoffline}), the iteration of the average state $\bar{x}_{t}$ satisfies
\begin{align}\label{iteration}
    \bar{x}_{t+1}=\bar{x}_{t} -\frac{\eta_t}{N}\sum_{i=1}^{N} \widetilde{\nabla}f_{i,t}(x_{i,t}).
\end{align}
Summing over $i\in [N]$, $\eta_t( f(\bar{x}_t)-f(x^*))$ can be bounded as
\begin{align*}
    &\eta_t (f(\bar{x}_t)-f(x^*)) \le N \langle\bar{x}_{t}-\bar{x}_{t+1}, \bar{x}_{t}-x^{*}\rangle\\
    & + \sum_{i=1}^{N} \eta_t \langle\theta_{i, t}, \bar{x}_{t}-x^{*}\rangle +\frac{L\eta_t}{2} \| \bar{x}_{t}-x_{i,t} \|^{2}.
\end{align*}
Note that 
\begin{align*}
    &\langle\bar{x}_{t}-\bar{x}_{t+1}, \bar{x}_{t}-x^{*}\rangle=\frac{1}{2}(\|\bar{x}_{t+1}-\bar{x}_{t}\|^2\\
    &+\|\bar{x}_{t}-x^{*}\|^2-\|\bar{x}_{t+1}-x^{*}\|^2).
\end{align*}
By (\ref{iteration}) and the boundedness of $\widetilde{\nabla}f_{i}(x_{i,.t})$, the validity of the lemma is ensured.
\subsection{Proof of Lemma \ref{lem_statebound}}\label{proofstate}
\textit{Proof:} By (\ref{boundofaverage}), we have that
\begin{align*}
    &z_n\|\bar{x}_{n+1}-x^*\|^{2}\le z_1\|\bar{x}_{1}-x^*\|^{2}
    +\sum_{t=1}^{n}z_t\eta_t(f(\bar{x}_t)-f(x^*))\\
    &\le z_1\|\bar{x}_{1}-x^*\|^{2} + \sum_{t=1}^{n}z_1\eta_t^2\lambda_t^2+\frac{2}{N} \sum_{i=1}^{N}Lz_1\eta_t \| \bar{x}_{t}-x_{i,t} \|^{2}\\
    &+\frac{1}{N} \sum_{i=1}^{N} 2z_t\eta_t \langle\theta_{i, t}, \bar{x}_{t}-x^{*}\rangle 
\end{align*}
where the inequality follows from the decreasing property of $z_t$ and 
\begin{align*}
    &\sum_{t=1}^{n} z_t\|\bar{x}_{t}-x^*\|^{2}- z_t\|\bar{x}_{t+1}-x^*\|^{2}\\
    &\le z_1\|\bar{x}_{1}-x^*\|^{2}- z_1\|\bar{x}_{2}-x^*\|^{2}\\
    &+\sum_{t=2}^{n} z_{t-1}\|\bar{x}_{t}-x^*\|^{2}- z_t\|\bar{x}_{t+1}-x^*\|^{2}\\
    & =z_1\|\bar{x}_{1}-x^*\|^{2}-z_n\|\bar{x}_{n+1}-x^*\|^{2}.
\end{align*}
Next we bound the square of the network error multiplying $\eta_t$. By Lemma \ref{lem_networkerr}, we obtain that
\textcolor{black}{\begin{equation}\label{networkerr^2}
    \begin{aligned}
    &\|x_{i,t}-\bar{x}_t\|^2\le  2(N\gamma)^{2}\Big(\Big.( \beta^{t-1} R_1)^2+
 (\sum_{l=1}^{t-1}\beta^{t-l}\lambda_{l}\eta_l)^2\Big)\Big. .
\end{aligned}
\end{equation}
}The inequality results from $(a_1+\cdots+a_n)^2\le n(a^2+\cdots+a_n^2)$. By the Cauchy-Schwartz inequality, we can derive that
\textcolor{black}{
\begin{equation}\label{networkerr^2witheta}
    \begin{aligned}
        &\sum_{t=1}^{T}\eta_t(\sum_{l=1}^{t-1}\beta^{t-l}\lambda_{l}\eta_l)^2\le \sum_{t=1}^{T}\eta_t(\sum_{l=1}^{t-1}\beta^{t-l})\sum_{l=1}^{t-1}\beta^{t-l}(\lambda_{l}\eta_l)^2\\
    % &\le (\sum_{l=1}^{t-1}\beta^{t-l})\sum_{t=1}^{T}\eta_t\sum_{l=1}^{t-1}\beta^{t-l}(\lambda_{l}\eta_l)^2
    % \\
    &\le \frac{1}{1-\beta}\sum_{l=1}^{T-1}\sum_{t=l+1}^{T-1}\beta^{t-l}(\lambda_{l}\eta_l)^2\eta_t\le (\frac{1}{1-\beta})^2\sum_{l=1}^{T-1}\lambda_{l}^2\eta_l^3.
    \end{aligned}
\end{equation}
}Combining all bounds of $\eta_{t}\|\bar{x}_{t}-x^*\|^2$ paves a way to bound the term $\|\bar{x}_{n+1}-x^*\|^2$ by $\|\bar{x}_{1}-x^*\|^2$ inductively. This leads to
 the validity of the result.
\subsection{Proof of Lemma \ref{lem_stateboundbynoise}}\label{proofnoise}
First we bound the term $\sum_{t=1}^{n} 2\eta_t z_t\langle \theta_{i,t}^u,\bar{x}_t-x^{*}\rangle$. It is worth noting that the  $\{z_t\eta_t\langle \theta_{i,t}^u,\bar{x}_t-x^{*}\rangle\}_{t\in[T]}$ is a sequence of martingale differences since $x_{i,t+1}$ is determined under $\mathcal{F}_{t}$. To bound it, we introduce a Lemma.
\begin{lemma}[\cite{nguyen2023improved}]\label{lem_exponential}
    Let  $X$  be a random variable such that  $\mathbb{E}[X]=0$ and  $|X| \leq R$  almost surely. Then for  $0 \leq \lambda \leq \frac{1}{R}$
    \begin{align*}
        \mathbb{E}[\exp (\lambda X)] \leq \exp \left(\frac{3}{4} \lambda^{2} \mathbb{E}\left[X^{2}\right]\right).
    \end{align*}
\end{lemma}
Now we construct a exponential martingale, and bound it via the Markov's inequality. 
\begin{align*}
    S_n= \sum_{t=1}^{n} 2\eta_t z_t\langle \theta_{i,t}^u,\bar{x}_t-x^{*}\rangle-\frac{3}{16}\eta_t^2\mathbb{E}[\|\theta_{i,t}^u\|^2|\mathcal{F}_t].
\end{align*}
Note that $\eta_t\|\theta_{i,t}^u\|\le 2$ almost surely, Lemma \ref{lem_exponential} indicates that
\begin{align*}
   \mathbb{E}[\exp (S_n)|\mathcal{F}_{n}]\le \exp(S_{n-1}).
\end{align*}
Thus $\{S_n\}_{n\in [T]}$ is an exponential supermartingale. By the Markov's inequality, we have, for $n\in [T]$
\begin{align*}
    P\{S_n\ge \log \frac{1}{\delta}\}\le \delta \mathbb{E}[\exp (S_1)]\le \delta.
\end{align*}
Then, with probability at least $1-\delta$, for any $n\in [T]$,
\begin{align}\label{term_residue}
    &\sum_{t=1}^{n} 2\eta_t z_t\langle \theta_{i,t}^u,\bar{x}_t-x^{*}\rangle \nonumber \\
    &\le \sum_{t=1}^{n} \frac{3}{16}\eta_t^2\mathbb{E}[\|\theta_{i,t}^u\|^2|\mathcal{F}_t] +\log \frac{1}{\delta}.
\end{align}

\textit{Proof:} Before using Lemma \ref{lem_heavytaied}, we need to bound $\Delta_{n}$.  Let $E_n$ be the event that for all $k\in [n]$,
\begin{align*}
    \sum_{t=1}^{k-1} 2z_t\eta_t\langle \theta_{i,t},\bar{x}_{ t}-x^{\star}\rangle \le \frac{5\Delta_{1}}{4}+\log \frac{1}{\delta}.
\end{align*}
We now prove inductively that under event $E_n$,  $\Delta_{n}\le (9\Delta_1+5a_0)$. For $k = 1$, the
event happens with probability 1, and $\Delta_1$ is already bounded. Suppose that $\Delta_{k}\le (9\Delta_1+5a_0)$ is valid under event $E_k$ for all $k\le n$, we obtain
\begin{align*}
    &\|\bar{x}_{n+1}-x^*\|^2\le \frac{1}{z_n} (
    \frac{9}{4}\Delta_1+a_0)\\
    &= (4(9\Delta_1+5a_0)+a_0)(
    \frac{9}{4}\Delta_1+a_0)\\
    &\le (9\Delta_1+5a_0)^2
\end{align*}
where $a_0=\max \{1,\log \frac{1}{\delta}\}$. Thus,  $\Delta_{n}$ is bounded. 
By condition (\ref{condition}), it follows that for any $t\in [n+1]$
\begin{align*}
    &\|\nabla f_{i}(x_{i,t})\|\\
    &\le\|\nabla f_{i}(x_{i,t})-\nabla f_{i}(x^*)\|+\|\nabla f_{i}(x^*)\|\\
    &\le L\|x_{i,t}-x^*\|+\|\nabla f_{i}(x^*)\|\\
    &\le L\|x_{i,t}-\bar{x}_t\|+L\|\bar{x}_{t}-x^*\|+\|\nabla f_{i}(x^*)\|\\
    &\le L\Delta_{t}+LD+B^*\\
    &\le L(9\Delta_1+5a_0)+LD+B^*\le \frac{\lambda_t}{2}
\end{align*}
\textcolor{black}{where $D=N\gamma R_1+
 N\gamma\frac{\lambda}{(1-\beta)m}$} is the bound of the network error. The result above is enough to obtain the probability of $E_{n+1}$.
First we bound $\sum_{t=1}^{n} z_t\eta_t\langle \theta_{i,t}^b,\bar{x}_t-x^{*}\rangle$, by Lemma \ref{lem_heavytaied}, it yields that
\begin{align*}
    &\sum_{t=1}^{n} z_t\eta_t\langle \theta_{i,t}^b,\bar{x}_t-x^{*}\rangle\\
    &\sum_{t=1}^{n} z_t\eta_t\|\theta_{i,t}^b\|\|\bar{x}_t-x^{*}\|\\
    &\le\frac{1}{4} \sum_{t=1}^{n} \eta_t \lambda_t (\frac{\sigma}{\lambda_t})^p\\
    &\le \textcolor{black}{\frac{\sigma^p}{m\lambda^{p-1}} \sum_{t=1}^{n} \frac{1}{t^{\kappa+(p-1)\alpha}(1+\log t)^2}}.
\end{align*}
By (\ref{condition}), we have
 $\sum_{t=1}^{n} z_t\eta_t\langle \theta_{i,t}^b,\bar{x}_t-x^{*}\rangle \le \Delta_{1}^2$. 
 
 Next we bound the term mentioned in (\ref{term_residue}), using Lemma \ref{lem_heavytaied} again, we derive that
 \begin{align*}
     &\sum_{t=1}^{n}\eta_t^2\mathbb{E}[\|\theta_{i,t}^u\|^2|\mathcal{F}_t] \\&\le \textcolor{black}{\frac{16}{m^2}\lambda^{2-p}\sigma^{p}\sum_{t=1}^{n} \frac{1}{t^{2\kappa-2\alpha}(1+\log t)^2}}.
 \end{align*}
 By (\ref{condition}) and combining the results above, the validity of the Lemma is ensured inductively.
\subsection{Proof of Theorem \ref{thm_highprobconvergence}}\label{proof_thm}
\textit{Proof:} By the analysis in Lemmas \ref{lem_boundofaverage} and \ref{lem_stateboundbynoise}, condition (\ref{condition}) and the decreasing property of $z_t$, with probability at least $1-\delta$, we have
\textcolor{black}{
\begin{align*}
     &\sum_{t=1}^{T}\Big(\Big. \frac{N\eta_t^2\lambda_t^2}{2}+\sum_{i=1}^{N} \eta_t \langle\theta_{i, t}, \bar{x}_{t}-x^{*}\rangle \\
     &+\frac{L\eta_t}{2} \| \bar{x}_{t}-x_{i,t} \|^{2}\Big)\Big.\\
     &\le 2L(N\gamma R_1)^2\frac{1}{1-\beta}+2L(N\gamma )^2(\frac{1}{1-\beta})^2\sum_{t=1}^{T-1}\lambda_{t}^2\eta_t^3 \\
     &+\sum_{t=1}^{T}\Big( \frac{N\eta_t^2\lambda_t^2}{2}+  \sum_{i=1}^{N} \eta_t \langle\theta_{i, t}, \bar{x}_{t}-x^{*}\rangle \Big)\\
     &\le 4N\Delta_{1}^2+ N\log \frac{1}{\delta}
\end{align*}
where the last line of the inequality above holds due to the second line of (\ref{condition}).}
Summing (\ref{boundofaverage}) over $t\in[T]$ and using the decreasing property of the step size $\eta_t$, we derive the result.

\end{document}